\documentclass[12pt]{amsart}
\usepackage{multicol}
\usepackage{amssymb}
\newtheorem{thm}{Theorem}[section]
\newtheorem{lem}[thm]{Lemma}
\newtheorem{cor}[thm]{Corollary}
\newtheorem{prop}[thm]{Proposition}
\newtheorem{ex}[thm]{Example}

\newtheorem*{prob*}{Open problem}

\theoremstyle{definition}

\newtheorem{defi}[thm]{Definition}

\theoremstyle{remark}

\newtheorem{rem}[thm]{Remark}
\newtheorem*{rem*}{Remark}

\newcommand{\kringel}{\mathbin{\raise1pt\hbox{$\scriptstyle\circ$}}}
\newcommand{\pkt}{\mathbin{\raise0pt\hbox{$\scriptstyle\bullet$}}}

\newcommand{\ad}{\mathop{\rm ad}}

\newcommand{\Lf}{\mathfrak{f}}
\newcommand{\Lg}{\mathfrak{g}}

\newcommand{\Ln}{\mathfrak{n}}

\newcommand{\Lt}{\mathfrak{t}}

\newcommand{\ga}{\gamma}

\renewcommand{\phi}{\varphi}
\begin{document}


\title[Novikov structures]{Novikov algebras and Novikov structures on Lie algebras}

\author[D. Burde]{Dietrich Burde}
\author[K. Dekimpe]{Karel Dekimpe}
\author[K. Vercammen]{Kim Vercammen}
\address{Fakult\"at f\"ur Mathematik\\
Universit\"at Wien\\
  Nordbergstr. 15\\
  1090 Wien \\
  Austria}
\email{dietrich.burde@univie.ac.at}
\address{Katholieke Universiteit Leuven\\
Campus Kortrijk\\
8500 Kortrijk\\
Belgium}
\date{\today}
\email{Karel.Dekimpe@kuleuven-kortrijk.be}
\email{Kim.Vercammen@kuleuven-kortrijk.be}

\begin{abstract}
We study ideals of Novikov algebras and Novikov structures on
finite-dimensional Lie algebras. We present the first example of a
three-step nilpotent Lie algebra which does not admit a Novikov structure.
On the other hand we show that any free three-step nilpotent Lie algebra
admits a Novikov structure. We study the existence question also for Lie algebras
of triangular matrices.
Finally we show that there are families of Lie algebras of arbitrary high solvability 
class which admit Novikov structures.
\end{abstract}

\maketitle

\section{Introduction}

Let $k$ be a field of characteristic zero. A Novikov algebra and, more
generally, an LSA is defined as follows:

\begin{defi}
An algebra $(A,\cdot)$ over $k$ with product $(x,y) \mapsto x\cdot y$
is called a {\it left-symmetric algebra (LSA)}, if the product is
left-symmetric, i.e., if the identity
\begin{equation}\label{lsa1}
x\cdot (y\cdot z)-(x\cdot y)\cdot z= y\cdot (x\cdot z)-(y\cdot x)\cdot z
\end{equation}
is satisfied for all $x,y,z \in A$. The algebra is called {\it Novikov}, if
in addition
\begin{equation}\label{nov}
(x\cdot y)\cdot z=(x\cdot z)\cdot y
\end{equation}
is satisfied.
\end{defi}

Denote by $L(x), R(x)$ the left, respectively right multiplication operator in
the algebra  $(A,\cdot)$. Then an LSA is a Novikov algebra if the right
multiplications commute:
\begin{align*}
[R(x),R(y)] & = 0.
\end{align*}

It is well known that LSAs are Lie-admissible algebras: the commutator
\begin{equation}\label{lsa2}
[x,y]=x\cdot y -y\cdot x
\end{equation}
defines a Lie bracket. The associated Lie algebra is denoted by $\Lg_A$.
The adjoint operator can be expressed by  $\ad(x)=L(x)-R(x)$. \\
If $A$ is a Novikov algebra, then we obtain, by expanding
the condition $0=[R(x),R(y)]=[L(x)-\ad(x),L(y)-\ad(y)]$, the following operator 
identity:

\begin{equation}\label{linearequation}
L([x,y])+\ad ([x,y])-[\ad(x),L(y)]-[L(x),\ad (y)] =0.
\end{equation}

\begin{defi}\label{affine}
An {\it affine structure} on a Lie algebra
$\Lg$ over $k$ is a left-symmetric product $\Lg \times \Lg \rightarrow \Lg$
satisfying \eqref{lsa2} for all $x,y \in \Lg$.
If the product is Novikov, we say that $\Lg$ admits a {\it Novikov structure}.
\end{defi}

A given Lie algebra need not admit a Novikov structure, or an affine structure.
The existence question for affine structures is very hard in general.
It is more accessible for Novikov structures. For results, background and references see,
for example, \cite{BEN},\cite{BUD}.

\section{Ideals in Novikov algebras}

In this section we will present some structure theory concerning ideals in
Novikov algebras. For related results in this direction see also \cite{BN}, \cite{ZEL}. 
We start with two identities which are similar to the Jacobi identity for Lie algebras. 

\begin{lem}\label{2.3}
Let $(A,\cdot)$ be a Novikov algebra. Then we have, for
all $x,y,z \in A$:
\begin{align*}
[x,y]\cdot z + [y,z]\cdot x + [z,x]\cdot y & = 0, \\
x \cdot [y,z] + y\cdot [z,x] + z\cdot [x,y] & = 0. \\
\end{align*}
\end{lem}
For a proof, see \cite{BU1}.
Next we show that the product of two ideals is again an ideal.

\begin{lem}
Let $(A,\cdot)$ be a Novikov algebra and $I,J$ be two-sided ideals of $A$.
Then $I \cdot J$ is also a two-sided ideal of $A$.
\end{lem}

\begin{proof}

Let $a \in A$, $x \in I$ and $y \in J$. Then the identity
\[
a \cdot (x \cdot y)=(a \cdot x) \cdot y + x \cdot (a \cdot y)-(x \cdot a)\cdot y
\]
shows that $a \cdot (x \cdot y) \in I \cdot J$. Because of
$(x \cdot y)\cdot a = (x \cdot a)\cdot y$ we also have
$(x \cdot y)\cdot a\in I \cdot J$.
\end{proof}

We also show that the commutator of two ideals is again an ideal.

\begin{lem}
Let $(A, \cdot)$ be a Novikov algebra and assume that $I,J$ are
two-sided ideals of $A$. Then $[I,J]$ is also a two-sided ideal of $A$.
\end{lem}

\begin{proof}
Let $a \in A$, $x \in I$ and $y \in J$. The operator identity (\ref{linearequation})
implies that
\begin{align*}
0 & = [x,y] \cdot a + [[x,y],a] - [x,y \cdot a] + y \cdot [x,a] - x \cdot
[y,a] + [y,x \cdot a]\\
  & = [x,y] \cdot a + [[x,y],a] - [x,y \cdot a] + [y,x \cdot a]+
(y \cdot [x,a] + x \cdot [a,y] + a\cdot [y,x] )- a \cdot [y,x]\\
  & = [x,y] \cdot a + [[x,y],a] - [x,y \cdot a] + [y,x \cdot a]+ a \cdot [x,y]\\
  & = [x,y] \cdot a + [[x,y],a] - [x,y \cdot a] + [y,x \cdot a]+ [x,y]\cdot
a + [a,[x,y]]\\
  & = 2 [x,y] \cdot a+ [y,x \cdot a]  - [x,y \cdot a].
\end{align*}
The term in brackets above vanishes because of lemma $\ref{2.3}$.
From this we deduce
\begin{align*}
[x,y]\cdot a & = \frac{1}{2} ([x,y \cdot a]-[y,x \cdot a]) \in [I,J], \\
a\cdot [x,y] & = [x,y] \cdot a  + [a,[x,y]]\in [I,J],
\end{align*}
which was to be shown.
\end{proof}

Let $(A,\cdot)$ be a Novikov algebra. Denote by
\begin{align*}
\gamma_1(A)& =\gamma_1(\Lg_A)=A \\
\gamma_{i+1}(A)& =\gamma_{i+1}(\Lg_A)=[A,\gamma_{i}(A)]
\end{align*}
the terms of the lower central series of $A$, respectively $\Lg_A$.
Furthermore denote by
\begin{align*}
A^{(0)}& =\Lg_A^{(0)}=A \\
A^{(i+1)}& =\Lg_A^{(i+1)}=[A^{(i)},A^{(i)}]
\end{align*}
the terms of the derived series of $A$, respectively $\Lg_A$. Then the above
lemma immediately implies the following result:

\begin{cor}\label{lowercentral}
Let $(A,\cdot)$ be a Novikov algebra. Then all $\gamma_i(A)$, and all $A^{(i)}$ are 
two-sided ideals of $A$.
\end{cor}

The ideals of the lower central series satisfy the following property.

\begin{lem}
Let $(A,\cdot)$ be a Novikov algebra. Then we have 
\[
\gamma_{i+1}(A) \cdot \gamma_{j+1}(A) \subseteq \gamma_{i+j+1}(A) 
\]
for all $i,j\ge 0$.
\end{lem}
\begin{proof}
We will show this by induction on $i\ge 0$.
The case $i=0$ follows from the fact that $\gamma_{j+1}(A)$ is an ideal in $A$, 
see corollary $\ref{lowercentral}$. Assume now that 
$\ga_k(A)\cdot \ga_{j+1}(A)\subseteq \ga_{k+j}(A)$ for all $k=1,\ldots ,i $.

Let $x\in \ga_1(A)$, $y\in \ga_i(A)$ and $z\in \ga_{j+1}(A)$. We have to show
that $[x,y]\cdot z \in \ga_{i+j+1}(A)$. 
The first identity of lemma $\ref{2.3}$
says, that $[x,y]\cdot z + [y,z]\cdot x + [z,x]\cdot y  = 0$. By $\eqref{lsa2}$ we have
\[
[z,x]\cdot y = y\cdot [z,x]+[[z,x],y].
\]
Then $y\cdot [z,x] \in \ga_{i+j+1}(A)$ by induction hypothesis, and
$[[z,x],y] \in \ga_{i+j+2}(A)$. It follows that $[z,x]\cdot y \in  \ga_{i+j+1}(A)$.
Similarly we obtain $[y,z]\cdot x \in  \ga_{i+j+1}(A)$. Now the first identity of
lemma $\ref{2.3}$ implies $[x,y]\cdot z\in  \ga_{i+j+1}(A)$.
\end{proof}

Denote the center of a Novikov algebra $A$ by $Z(A)=\{ x\in A\mid x\cdot y
=y\cdot x \mbox{ for all } y \in A\}$. Note that $Z(A)$ is also the center of the 
associated Lie algebra $\Lg_A$.

\begin{lem}\label{2.4}
Let $(A, \cdot)$ be a Novikov algebra. Then $Z(A) \cdot [A,A] =[A,A]\cdot Z(A)=0$.
\end{lem}

\begin{proof}
Let $a,b \in A$ and $z \in Z(A)$. Again by lemma $\ref{2.3}$
we have
\[
z \cdot [a,b] + a \cdot [b,z] + b \cdot [z,a] = 0.
\]
Since $z$ is also in the center of the associated Lie algebra of $A$
we obtain $z \cdot [a,b]=0$. Furthermore we have
\begin{align*}
0 & = [z,[b,a]] = z\cdot [b,a] - [b,a]\cdot z \\
  & = [a,b]\cdot z.
\end{align*}
\end{proof}

\begin{lem}\label{2.5}
Let $(A, \cdot)$ be a Novikov algebra. Then $Z(A)$ is a two-sided ideal
of $A$.
\end{lem}

\begin{proof}
Let $z \in Z(A)$. For any $b \in A$ we have the following
two identities
\begin{align*}
[L(b),L(z)] &= L([b,z])= 0,\\
[R(b),R(z)] &= 0.
\end{align*}
Because $z \in Z(A)$ we have $R(z) = L(z)$. Hence we also have
$[L(b),R(z)]= 0$, so that
\begin{align*}
0 & = [L(b)-R(b),R(z)] = [\ad(b), R(z)].
\end{align*}
In particular it follows $[b, a\cdot z] - [b,a]\cdot z =0$ for all $a\in A$.
By Lemma $\ref{2.4}$ we have $[b,a]\cdot z = 0$, hence
$[b, a\cdot z]=0$. Because this is true for every $b \in A $, we can
conclude that $a \cdot z \in Z(A)$.
Since $z \in Z(A)$ we also have $z \cdot a \in Z(A)$.
\end{proof}

Let $Z_1(A)=Z(A)$ and define $Z_{i+1}(A)$ by the identity 
$Z_{i+1}(A)/Z_{i}(A)= Z(A/Z_i(A))$. Note that the $Z_i(A)$ are the terms of the 
upper central series of the associated Lie algebra $\Lg_A$.
As an immediate consequence of the previous lemma, we obtain

\begin{cor}
Let $(A,\cdot)$ be a Novikov algebra. Then all terms $Z_i(A)$ of the upper central series
of $A$ are two-sided ideals of $A$.
\end{cor}

Denote by $(x,y,z)=x\cdot (y\cdot z)-(x\cdot y)\cdot z$ the associator of three
elements in $A$.

\begin{lem}
Let $A$ be a Novikov algebra and one of the elements $x,y,z$ in $Z(A)$. Then
$(x,y,z)=0$.
\end{lem}

\begin{proof}
In any LSA we have the identity
\[
(x,y,z)=x\cdot [y,z]+[z,x\cdot y]+[x,z]\cdot y.
\]
If $z\in Z(A)$, then this implies $(x,y,z)=0$. If $y\in Z(A)$ then
also $x\cdot y \in Z(A)$ by lemma $\ref{2.5}$, and $[A,A]\cdot Z(A)=0$ by lemma
$\ref{2.4}$. Hence the above identity implies  $(x,y,z)=0$. The same argument
shows the claim for $x\in Z(A)$.
\end{proof}

\section{Novikov structures on 3-step nilpotent Lie algebras}

In \cite[Remark 4.11]{BUD} it was questioned whether or not there exists a 3-step
nilpotent Lie algebra not admitting a Novikov structure. In the same paper it was 
shown that
a Novikov structure does exist when the 3-step nilpotent Lie algebra $\Lg$ can be generated
by at most 3 elements. This result was obtained by first considering a Novikov structure
on the free 3-step nilpotent Lie algebra $\Lf$ on 3 generators and then it was shown that
$\Lg$ could be realized as a quotient $\Lg=\Lf/I$, where $I$ is an ideal of $\Lf$ seen as
a Novikov algebra.

\medskip

Having this in mind, we first study the free 3-step nilpotent case.
\begin{prop}
Let $\Lg$ be a free 3-step nilpotent Lie algebra on $n$ generators
$x_1,x_2,\ldots,x_n$. Then $\Lg$ admits a Novikov structure.
\end{prop}

\begin{proof}
As a vector space, $\Lg$ has a basis
\[ x_1,x_2,\ldots, x_n, \]
\[ y_{i,j}=[x_i,x_j],\;(1\leq i < j \leq n),\;\]
\[z_{i,j,k}=[x_i, y_{j,k}],\;(1\leq j < k \leq n, \; 1\leq  i \leq k
\leq n).\]
Note that in case $i>k$, we have that
\[
z_{i,j,k}=[x_i,y_{j,k}] = [x_i,[x_j,x_k]]= -[x_j,[x_k,x_i]]-[x_k,[x_i,x_j]]=
-z_{j,k,i}+z_{k,j,i}.
\]
A Novikov structure on $\Lg$ is defined by \\[0.2cm]
\begin{itemize}
\item If $n\geq i > j \geq 1$ then $x_i\cdot x_j=-y_{j,i}$. \\
\item If $1\leq i \leq  j < k\leq n $ then $x_i\cdot y_{j,k}=
\frac{z_{i,j,k}}{2}$. \\
If $1\leq j  < i < k\leq n $ then $x_i\cdot y_{j,k}=
-\frac{z_{j,i,k}}{2}+z_{i,j,k}$. \\
If $1\leq j  < k \leq i  \leq n $ then $x_i\cdot y_{j,k}= z_{i,j,k}$. \\
\item If $1\leq i \leq j < k \leq n$, then $y_{j,k}\cdot
x_i=-\frac{z_{i,j,k}}{2}$.\\
If $1\leq j < i < k \leq n$, then $y_{j,k}\cdot x_i=-\frac{z_{j,i,k}}{2}$.\\
\end{itemize}
All other products are zero. \\[0.5cm]
By considering each of the above cases it is easy to see that the identity
$[a,b] = a \cdot b - b \cdot a$ holds for all basis elements $a$ and $b$.
We have to show the following two other identities:
\begin{eqnarray*}(a \cdot b) \cdot c - a \cdot (b \cdot c) - (b \cdot a)
\cdot c + b \cdot (a \cdot c) &=& 0,  \\
(a \cdot b)\cdot c - (a \cdot c)\cdot b &=& 0,
\end{eqnarray*}
for all basis elements $a,b$ and $c$. It is clear that 
we only have to consider the case where $a= x_i$, $b=x_j$ and $c=x_k$: otherwise
the two identities will be trivially satisfied, because any product of
the form $y_{i,j} \cdot y_{k,l}$ is zero, and any product that involves
an element $z_{i,j,k}$ is also zero.
For the first condition we will consider the case $1 \leq k<i<j \leq n$. We have
\begin{align*}
& (x_i \cdot x_j) \cdot x_k - x_i \cdot (x_j \cdot x_k) - (x_j \cdot
x_i) \cdot x_k + x_j \cdot (x_i \cdot x_k)\\
& = - x_i \cdot (-y_{k,j}) - (-y_{i,j})\cdot x_k + x_j \cdot (-y_{k,i})\\
& = - \frac{z_{k,i,j}}{2} + z_{i,k,j}- \frac{z_{k,i,j}}{2}-z_{j,k,i}\\
& = - z_{k,i,j} + z_{i,k,j}+z_{k,i,j}-z_{i,k,j}\\
&= 0.
\end{align*}

Similarly the other cases can be
treated.
For the second condition we consider the case $1 \leq j<k<i \leq n$. We have
\begin{align*}
(x_i \cdot x_j)\cdot x_k - (x_i \cdot x_k)\cdot x_j & = -y_{j,i} \cdot x_k + y_{k,i} 
\cdot x_j\\
& = \frac{z_{j,k,i}}{2} - \frac{z_{j,k,i}}{2}\\
& =0.
\end{align*}
Similarly the other cases can be shown. It follows that the product defines a Novikov 
structure on $\Lg$.
\end{proof}
As a motivation for what follows, we provide a detailed description in the four 
generator case.

\begin{ex}\label{to-use}
Let $n=4$. Then $\dim \Lg=30$.
The nonzero Lie brackets and Novikov products are given as follows.
\end{ex}

\begin{multicols}{3}
$[x_{1},x_{2}]=y_{1,2}$\\
$[x_{1},x_{3}]=y_{1,3}$\\
$[x_{1},x_{4}]=y_{1,4}$\\
$[x_{1},y_{1,2}]=z_{1,1,2}$\\
$[x_{1},y_{1,3}]=z_{1,1,3}$\\
$[x_{1},y_{1,4}]=z_{1,1,4}$\\
$[x_{1},y_{2,3}]=z_{1,2,3}$\\
$[x_{1},y_{2,4}]=z_{1,2,4}$\\
$[x_{1},y_{3,4}]=z_{1,3,4}$\\
$[x_{2},x_{3}]=y_{2,3}$\\
$[x_{2},x_{4}]=y_{2,4}$\\
$[x_{2},y_{1,2}]=z_{2,1,2}$\\
$[x_{2},y_{1,3}]=z_{2,1,3}$\\
$[x_{2},y_{1,4}]=z_{2,1,4}$\\
$[x_{2},y_{2,3}]=z_{2,2,3}$\\
$[x_{2},y_{2,4}]=z_{2,2,4}$\\
$[x_{2},y_{3,4}]=z_{2,3,4}$\\
$[x_{3},x_{4}]=y_{3,4}$\\
$[x_{3},y_{1,2}]=-z_{1,2,3} + z_{2,1,3}$\\
$[x_{3},y_{1,3}]=z_{3,1,3}$\\
$[x_{3},y_{1,4}]=z_{3,1,4}$\\
$[x_{3},y_{2,3}]=z_{3,2,3}$\\
$[x_{3},y_{2,4}]=z_{3,2,4}$\\
$[x_{3},y_{3,4}]=z_{3,3,4}$\\
$[x_{4},y_{1,2}]=-z_{1,2,4} + z_{2,1,4}$\\
$[x_{4},y_{1,3}]=-z_{1,3,4} + z_{3,1,4}$\\
$[x_{4},y_{1,4}]=z_{4,1,4}$\\
$[x_{4},y_{2,3}]=-z_{2,3,4} + z_{3,2,4}$\\
$[x_{4},y_{2,4}]=z_{4,2,4}$\\
$[x_{4},y_{3,4}]=z_{4,3,4}$
\end{multicols}

\begin{multicols}{3}
$x_{1}\cdot y_{1,2}=\frac{z_{1,1,2}}{2}$\\
$x_{1}\cdot y_{1,3}=\frac{z_{1,1,3}}{2}$\\
$x_{1}\cdot y_{1,4}=\frac{z_{1,1,4}}{2}$\\
$x_{1}\cdot y_{2,3}=\frac{z_{1,2,3}}{2}$\\
$x_{1}\cdot y_{2,4}=\frac{z_{1,2,4}}{2}$\\
$x_{1}\cdot y_{3,4}=\frac{z_{1,3,4}}{2}$\\
$x_{2}\cdot x_{1}=-y_{1,2}$\\
$x_{2}\cdot y_{1,2}=z_{2,1,2}$\\
$x_{2}\cdot y_{1,3}=\frac{-z_{1,2,3}}{2}+ z_{2,1,3}$\\
$x_{2}\cdot y_{1,4}=\frac{-z_{1,2,4}}{2} + z_{2,1,4}$\\
$x_{2}\cdot y_{2,3}=\frac{z_{2,2,3}}{2}$\\
$x_{2}\cdot y_{2,4}=\frac{z_{2,2,4}}{2}$\\
$x_{2}\cdot y_{3,4}=\frac{z_{2,3,4}}{2}$\\
$x_{3}\cdot x_{1}=-y_{1,3}$\\
$x_{3}\cdot x_{2}=-y_{2,3}$\\
$x_{3}\cdot y_{1,2}=-z_{1,2,3} + z_{2,1,3}$\\
$x_{3}\cdot y_{1,3}=z_{3,1,3}$\\
$x_{3}\cdot y_{1,4}=\frac{-z_{1,3,4}}{2}+ z_{3,1,4}$\\
$x_{3}\cdot y_{2,3}=z_{3,2,3}$\\
$x_{3}\cdot y_{2,4}=\frac{-z_{2,3,4}}{2}+ z_{3,2,4}$\\
$x_{3}\cdot y_{3,4}=\frac{z_{3,3,4}}{2}$\\
$x_{4}\cdot x_{1}=-y_{1,4}$\\
$x_{4}\cdot x_{2}=-y_{2,4}$\\
$x_{4}\cdot x_{3}=-y_{3,4}$\\
$x_{4}\cdot y_{1,2}=-z_{1,2,4} + z_{2,1,4}$\\
$x_{4}\cdot y_{1,3}=-z_{1,3,4} + z_{3,1,4}$\\
$x_{4}\cdot y_{1,4}=z_{4,1,4}$\\
$x_{4}\cdot y_{2,3}=-z_{2,3,4} + z_{3,2,4}$\\
$x_{4}\cdot y_{2,4}=z_{4,2,4}$\\
$x_{4}\cdot y_{3,4}=z_{4,3,4}$\\
$y_{1,2}\cdot x_{1}=-\frac{z_{1,1,2}}{2}$\\
$y_{1,3}\cdot x_{1}=-\frac{z_{1,1,3}}{2}$\\
$y_{1,3}\cdot x_{2}=-\frac{z_{1,2,3}}{2}$\\
$y_{1,4}\cdot x_{1}=-\frac{z_{1,1,4}}{2}$\\
$y_{1,4}\cdot x_{2}=-\frac{z_{1,2,4}}{2}$\\
$y_{1,4}\cdot x_{3}=-\frac{z_{1,3,4}}{2}$\\
$y_{2,3}\cdot x_{1}=-\frac{z_{1,2,3}}{2}$\\
$y_{2,3}\cdot x_{2}=-\frac{z_{2,2,3}}{2}$\\
$y_{2,4}\cdot x_{1}=-\frac{z_{1,2,4}}{2}$\\
$y_{2,4}\cdot x_{2}=-\frac{z_{2,2,4}}{2}$\\
$y_{2,4}\cdot x_{3}=-\frac{z_{2,3,4}}{2}$\\
$y_{3,4}\cdot x_{1}=-\frac{z_{1,3,4}}{2}$\\
$y_{3,4}\cdot x_{2}=-\frac{z_{2,3,4}}{2}$\\
$y_{3,4}\cdot x_{3}=-\frac{z_{3,3,4}}{2}$
\end{multicols}

When trying to find an example of a $3$-step nilpotent Lie algebra without a Novikov 
structure, we know from \cite{BUD} that such an example must have at least $4$ generators. 
Any $3$-step nilpotent Lie algebra on $4$ generators is a quotient $\Lg/I$ of the Lie 
algebra $\Lg$ described in example $\ref{to-use}$, where $I$ is an ideal of 
$\Lg$. So, in order to find such a Lie algebra $\Lg/I$ without
Novikov structure, we have to choose an $I$ which is certainly not an ideal of $\Lg$, 
seen as a Novikov algebra. The following proposition uses such an example.

\begin{prop}\label{counterexample}
Consider the following $3$-step nilpotent Lie algebra $\Lg$ on $4$ generators
of dimension $13$, with basis $(x_1,\ldots,x_{13})$ and non-trivial Lie brackets
\begin{align*}
[x_1,x_2] & = x_5, && [x_3,x_4]=-x_5, \\
[x_1,x_4] & = x_6,  &&[x_3,x_5]=-x_{11}, \\
[x_1,x_6] & = x_{10}, && [x_3,x_8]=x_{9}, \\
[x_1,x_7] & = x_{11},  &&[x_4,x_5]=-x_{12},\\
[x_1,x_8] & = x_{12}, && [x_4,x_6]=x_{9},\\
[x_2,x_3] & = x_7,&&  [x_4,x_7]=x_9+x_{13}.\\
[x_2,x_4] & = x_8, \\
[x_2,x_5] & = x_{13}, \\
[x_2,x_7]  &= x_{13}, \\
\end{align*}
This Lie algebra does not admit a Novikov structure.
\end{prop}

Note that $\Lg$ admits an affine structure since it is positively graded.

\begin{proof}
We will assume that $\Lg$ admits a Novikov structure and show that this leads to a
contradiction.
We express the adjoint operators $\ad(x_i)$ and the left (resp.\ right) multiplication 
operators $L(x_i) $ (resp.\ $R(x_i)$) as matrices with respect to the basis 
$x_1,x_2,\ldots, x_{13}$.
The adjoint operators $\ad(x_i)$ are given by the Lie bracktes of $\Lg$, while the 
left multiplication operators are unknown. We denote the $(j,k)$-th entry
of $L(x_i)$ by
\[ 
L(x_i)_{j,k} = x^i_{j,k}.
\]
We use the convention that the $j$-th column of $L(x_i)$ gives the coordinates of
$L(x_i)(x_j)$. Note that once the entries of the left multiplication operators are chosen,
the right multiplication operators are given by $R(x_i)_{j,k}=x^k_{j,i}$.
We have to satisfy all relations given by $\eqref{lsa1}$, $\eqref{nov}$ and
$\eqref{lsa2}$, where $x$, $y$ and $z$ run over all basis vectors. This leads to a 
huge system of quadratic equations in the variables $x^i_{j,k}$ for
$1\leq i,j,k \leq 13$, summing up to a total of $13^3=2197$ variables. We need to
show that these equations are contradictory. At first sight, this seems to be a rather 
hopeless task. 
However, we can use our knowledge on ideals in a Novikov algebra. Then we find that a 
lot of unknowns $x^{i}_{j,k}$ already have to be zero. In the table below, we list the 
triples $(i,j,k)$ for which we already know that $x^i_{j,k}=0$:
\begin{center}
\begin{tabular}{llll}
$1\leq i \leq 13$, & $1\leq  j \leq 4$, & $5 \leq k \leq 13$, & because 
$\gamma_2(\Lg)$ is an ideal.\\
$1\leq i \leq 13$, & $5\leq  j \leq 8$, & $9 \leq k \leq 13$, & because 
$\gamma_3(\Lg)$ is an ideal.\\
$5\leq i \leq 13$, & $1\leq  j \leq 4$, & $1 \leq k \leq 4$,  & because 
$\gamma_2(\Lg)$ is an ideal.\\
$5\leq i \leq 13$, & $5\leq  j \leq 8$, & $5 \leq k \leq 8$, &
                because $\gamma_2(\Lg)\cdot\gamma_2(\Lg)\subseteq \gamma_3(\Lg)$.\\
$5\leq i \leq 13$, & $9\leq  j \leq 13$, & $9 \leq k \leq 13$, &
                because $\gamma_2(\Lg)\cdot\gamma_3(\Lg)=0$.\\
$9\leq i \leq 13$, & $5\leq  j \leq 8$, & $1 \leq k \leq 4$,  & because 
$\gamma_3(\Lg)$ is an ideal.\\
$9\leq i \leq 13$, & $9\leq  j \leq 13$, & $5 \leq k \leq 8$, & because 
$\gamma_3(\Lg)\cdot\gamma_2(\Lg)=0$.
\end{tabular}
\end{center}
It follows that $1421$ of the $x^i_{j,k}$ have to be zero, leaving us with 
$776$ variables. \\
On the other hand the conditions  
\[ 
\ad(x_i)= L(x_i)-R(x_i),\;\;1\leq i \leq 13.
\]
yield a (large but very simple) system of linear equations; allowing us to 
determine $352$ variables $x^i_{j,k}$ in dependance of the remaining $776-352= 424$ ones.\\
To get a further reduction we use that  
\[ 
x_i\cdot [x_j,x_k]+ x_j\cdot [x_k,x_i]+ x_k \cdot [x_i,x_j]=0,\quad
1 \leq i < j < k \leq 13,
\]
which is the same as
\[ 
L(x_i) (\ad(x_j) x_k) +
L(x_j) (\ad(x_k) x_i) +
L(x_k) (\ad(x_i) x_j)=0, \quad
1 \leq i < j < k \leq 13.
\]
Again this leads to a system of linear equations, this time specifying $156$ unknowns in
terms of the other ones, leaving $424-156=268$ variables.\\
Now, we consider the operator identity $\eqref{linearequation}$, i.e.,
\[
L([x_i,x_j])+\ad ([x_i,x_j])-[\ad(x_i),L(x_j)]-[L(x_i),\ad (x_j)] =0,
\quad 1\leq i<j\leq 13.
\]
Note that for any pair $(i,j)$, we can write $[x_i,x_j]$ as a linear combination of
the $x_k$, $1\leq k \leq 13$. Hence we can also write $L([x_i,x_j])$ as the corresponding
linear combination of the $L(x_k)$. Doing this, we obtain another system of linear 
equations, determining $210$ extra variables, leaving $268-210=58$ free variables. \\
Finally we use that the right multiplications have to commute, i.e.,
\[ 
R(x_i) R(x_j)-R(x_j) R(x_i) =0, \;\; 1\leq i < j \leq 13.
\]
This yields a system of quadratic equations, which is immediately contradictory. In fact,
when taking $i=1$ and $j=2$, one obtains the equation $0=\frac18$, which is the 
desired contradiction.
\end{proof}

\section{The (non) existence of Novikov structures on triangular matrix algebras}

One of the most fundamental examples for solvable, resp.\ nilpotent Lie algebras
are the Lie algebras of upper-triangular, resp.\ strictly upper triangular
matrices of size $n$ over a field $k$, which we denote by $\Lt(n,k)$, resp.\ $\Ln(n,k)$.
It is therefore natural to ask, which of those Lie algebras admit a Novikov structure.  
It turns out that such structures exist only in very small dimensions. 

\begin{prop}
The Lie algebra $\Ln(n,k)$ admits a Novikov structure if and only if $n\leq 4$.
\end{prop}

\begin{proof}
If $n\leq 4$, the Lie algebra $\Ln(n,k)$ is abelian ($n=2$), $2$-step nilpotent ($n=3$) or
$3$-step nilpotent and generated by $3$ elements ($n=4$). In any of these cases, 
we know that a Novikov structure exists. \\
Now let $n>4$ and suppose that $\Ln(n,k)$ admits a Novikov structure. 
Denote by $e_{i,j}$ the elementary matrices, which have a $1$ on  the $(i,j)$-th position 
and a zero elsewhere. The $e_{i,j}$ with  $1\leq i < j  \leq n$ form a basis 
of $\Ln(n,k)$. The Lie bracket is given by
\[ 
[e_{i,j}, e_{k,l}]=\delta_{j,k}e_{i,l}-\delta_{i,l} e_{k,j}. 
\]
Assume that $(A,\cdot)$ defines a Novikov structure on $\Ln(n,k)$.
Then some easy calculations, using lemma $\ref{2.3}$ and identity $\eqref{lsa2}$ yield:
\begin{align*}
e_{1,2}\cdot [e_{3,4},e_{4,5}]+e_{3,4}\cdot [e_{4,5},e_{1,2}] + e_{4,5}
\cdot [e_{1,2},e_{3,4}] & =0 \quad \Rightarrow \; e_{1,2}\cdot e_{3,5}=0 \\
e_{3,4}\cdot [e_{4,5},e_{1,3}] + e_{4,5}\cdot [e_{1,3},e_{3,4}] +
e_{1,3}\cdot [e_{3,4},e_{4,5}] & = 0 \quad 
\Rightarrow \; e_{1,3}\cdot e_{3,5}= - e_{4,5}\cdot e_{1,4} \\
e_{3,4}\cdot [e_{4,5},e_{2,3}] + e_{4,5}\cdot [e_{2,3},e_{3,4}] +
e_{2,3}\cdot [e_{3,4},e_{4,5}] & = 0 \quad 
\Rightarrow \; e_{2,3}\cdot e_{3,5}= - e_{4,5}\cdot e_{2,4} \\
e_{1,2}\cdot [e_{4,5},e_{2,4}] + e_{4,5}\cdot [e_{2,4},e_{1,2}] +
e_{2,4}\cdot [e_{1,2},e_{4,5}] & = 0 \quad
\Rightarrow \; e_{1,2}\cdot e_{2,5}= - e_{4,5}\cdot e_{1,4}=e_{1,3}\cdot e_{3,5} \\
e_{1,2}\cdot [e_{2,3},e_{3,5}] + e_{2,3}\cdot [e_{3,5},e_{1,2}] +
e_{3,5}\cdot [e_{1,2},e_{2,3}] & = 0 \quad 
\Rightarrow \; e_{1,2}\cdot e_{2,5}= - e_{3,5}\cdot e_{1,3}=e_{1,3}\cdot e_{3,5} \\
[e_{1,3},e_{3,5}]- e_{1,3}\cdot e_{3,5} + e_{3,5}\cdot e_{1,3} & = 0\quad
\Rightarrow \; e_{1,3}\cdot e_{3,5} = e_{1,5}/2 = - e_{3,5} \cdot e_{1,3} \\
[e_{1,4},e_{4,5}]- e_{1,4}\cdot e_{4,5}+e_{4,5}\cdot e_{1,4} & = 0\quad   
\Rightarrow \; e_{1,4} \cdot e_{4,5}=e_{1,5}/2 =  - e_{4,5}\cdot e_{1,4}.
\end{align*}
Applying the operator identity $\eqref{linearequation}$ for $x=e_{1,2}$, $y=e_{2,3}$ to
$z=e_{3,5}$, and using the above computations, we find
\begin{eqnarray*}
0 & = & (L([e_{1,2},e_{2,3}]) + \ad[e_{1,2},e_{2,3}] -[L(e_{1,2}),\ad(e_{2,3})]-
         [\ad(e_{1,2}),L(e_{2,3})]) (e_{3,5})\\
  & = & e_{1,3}\cdot e_{3,5} + e_{1,5} -e_{1,2}\cdot e_{2,5} - [e_{1,2}, e_{2,3}\cdot e_{3,5}] \\
  & = & e_{1,5} + [e_{1,2}, e_{4,5}\cdot e_{2,4}]\\
  & = & e_{1,5} + [e_{1,2}, e_{2,4}\cdot e_{4,5} + [e_{4,5},e_{2,4}]].
\end{eqnarray*}
It follows that $[e_{1,2}, e_{2,4}\cdot e_{4,5} ]=0$. Furthermore it follows that
\begin{eqnarray*}
0 & = & (L([e_{1,2},e_{2,4}]) + \ad[e_{1,2},e_{2,4}] -[L(e_{1,2}),\ad(e_{2,4})]-
         [\ad(e_{1,2}),L(e_{2,4})]) (e_{4,5})\\
 & = & e_{1,4}\cdot e_{4,5} + e_{1,5} - e_{1,2}\cdot e_{2,5} +
 [e_{2,4}, e_{1,2}\cdot e_{4,5}] - [e_{1,2}, e_{2,4}\cdot e_{4,5} ]\\
 & =  & e_{1,5} + [e_{2,4}, e_{1,2} \cdot e_{4,5} ].
 \end{eqnarray*}
However, this is impossible, since $e_{1,5}\not \in [e_{2,4},\Ln(n,k)]$. This
contradiction shows that there is no Novikov structure on $\Ln(n,k)$ when $n\geq 5$.
\end{proof}

As a consequence we can easily prove the following result.

\begin{prop}
The Lie algebra $\Lt(n,k)$  admits a Novikov structure if and only if $n\leq2$.
\end{prop}

\begin{proof}
It is easy to construct a Novikov structure on $\Lt(1,k)\cong k$ and on $\Lt(2,k)$.
For $n=3$ and $n=4$ it is not difficult to see by direct calculations 
that $\Lt(n,k)$ does not admit a Novikov structure. For $n\ge 5$ we can use the above
proposition. Assume that  $\Lt(n, k)$ admits a Novikov structure for $n\ge 5$. 
Then also  $[\Lt(n,k),\Lt(n,k)]=\Ln(n,k)$ admits a Novikov structure, in contradiction
to our previous proposition.
\end{proof}

\section{Novikov structures on k-step solvable Lie algebras}

A natural question is, whether there are families of Lie algebras
of solvability class $k$, which admit Novikov structures for all $k\ge 1$.
The same question for nilpotency class has an easy answer. Here the standard filiform
nilpotent Lie algebras with basis $(e_1,\ldots ,e_n)$ and brackets $[e_1,e_i]=e_{i+1}$
for $i=2,\ldots ,n-1$ admit Novikov structures. Hence they
provide examples of nilpotency class $k=n-1$, see \cite{BUD}.
The following result shows that there are indeed filiform nilpotent Lie algebras of arbitrary
solvability class, which admit Novikov structures. Define for every $n\ge 3$ a filiform
Lie algebra $\Lf_{\frac{9}{10},n}$ of dimension $n$ by
\begin{align*}
[e_1,e_j] & = e_{j+1}, \quad 2\le j\le n-1,\\[0.2cm]
[e_i,e_j] & = \frac{6(j-i)}{j(j-1)\binom{j+i-2}{i-2}}e_{i+j},
\quad 2\le i \le j; \; i+j \leq n
\end{align*}
In particular we have
\begin{align*}
[e_2,e_j] & = \frac{6(j-2)}{j(j-1)} e_{j+2}, \quad 3\le j\le n-3, \\[0.2cm]
[e_j,e_{j+1}]& =\frac{6(j-1)!(j-2)!}{(2j-1)!}e_{2j+1},
\quad 2\le j\le (n-1)/2 .
\end{align*}
Then $[e_2,e_3]=e_5$, $[e_2,e_4]=e_6$, $[e_2,e_5]=\frac{9}{10}e_7$, etc. Similar Lie
algebras were studied in \cite{BEN}. \\
To verify the Jacobi identity introduce a new basis $(f_1,\ldots ,f_n)$ by
\begin{align*}
f_1 & = 6e_1,\\
f_j & = \frac{1}{(j-2)!}e_j, \quad 2\le j\le n.
\end{align*}
Then the new brackets are given by
\begin{align*}
[f_i,f_j] & = 6(j-i)f_{i+j}, \quad 1\le i \le j; \; i+j\le n.
\end{align*}
Here the Jacobi identity is obvious.
\begin{prop}
For each $n\ge 3$ the Lie algebra $\Lf_{\frac{9}{10},n}$ admits a complete Novikov structure.
It is given by the following multiplication:
\begin{align*}
e_1\cdot e_j & = e_{j+1},\quad 2\le j\le n-1,\\[0.2cm]
e_i\cdot e_j & = \frac{6}{j\binom{j+i-2}{i-2}}e_{i+j},
\quad 2\le i ,j \le n, \; i+j\le n.
\end{align*}
\end{prop}

\begin{rem}
Note that $\Lf_{\frac{9}{10},n}$ is $k$-step solvable if $2^{k}\leq n+1 < 2^{k+1}$. Indeed,
\begin{align*}
\Lg^{(0)}& =\Lg, \\
\Lg^{(1)}& =[\Lg,\Lg]=\langle e_3,\ldots , e_n \rangle, \\
\Lg^{(i)}& =[\Lg^{(i-1)},\Lg^{(i-1)}] =\langle e_{2^{i+1}-1},\ldots , e_n \rangle.
\end{align*}
Hence these algebras can have arbitrary high solvability class.
\end{rem}

\begin{proof}
In the new basis $(f_1,\ldots ,f_n)$ the Novikov product is given by
\begin{align*}
f_i\cdot f_j & = 6(j-1)f_{i+j}, \quad 1\le i , j\leq n,\; i+j\le n.
\end{align*}
Now it is easy to verify the required identities.
We have
(using the covention that $f_m=0$ when $m>n$)
\begin{align*}
f_i\cdot f_j-f_j\cdot f_i & = 6(j-1)f_{i+j}-6(i-1)f_{i+j} \\
 & = 6(j-i)f_{i+j}=[f_i,f_j], \quad i,j\ge 1,
\end{align*}
so that \eqref{lsa2} is satisfied. We have
\begin{align*}
(f_i\cdot f_j)\cdot f_k & = 36(j-1)(k-1)f_{i+j+k},\\
(f_i\cdot f_k)\cdot f_j & = 36(k-1)(j-1)f_{i+j+k},\\
\end{align*}
so that \eqref{nov} is satisfied. Finally,
\begin{align*}
f_i\cdot (f_j\cdot f_k) -(f_i\cdot f_j)\cdot f_k & = 36\cdot k(k-1) f_{i+j+k},\\
f_j\cdot (f_i\cdot f_k) -(f_j\cdot f_i)\cdot f_k & = 36\cdot k(k-1) f_{i+j+k},\\
\end{align*}
so that \eqref{lsa1} is satisfied.
\end{proof}


\begin{thebibliography}{99}


\bibitem{BN} A. A. Balinskii; S. P. Novikov: {\it Poisson brackets of
hydrodynamic type, Frobenius algebras and Lie algebras}.
Sov.\ Math.\ Dokl.\ \textbf{32} (1985), 228-231.

\bibitem{BEN} Y. Benoist, {\it Une nilvari\'et\'e non affine}. J.\ Differential
Geom. \textbf{41} (1995), 21-52.

\bibitem{BU1} D. Burde: {\it Classical r-matrices and Novikov algebras}.
Geom.\ Dedicata \textbf{122}, No. 1 (2006), 145-157.

\bibitem{BUD} D. Burde, K. Dekimpe: {\it Novikov structures on
solvable Lie algebras}. J.\ Geom.\ Phys.\ {\bf 56}, Nr. 9, 1837-1855 (2006).

\bibitem{GD}  I. M. Gel'fand;  I. Dorfman: {\it Hamiltonian operators
and algebraic structures related to them}. Funct.\ Anal.\ Appl.\
\textbf{13} (1980), 248-262.

\bibitem{ZEL} E. Zelmanov: {\it  On a class of local translation invariant
Lie algebras}. Soviet Math. Dokl. {\bf 35} (1987), 216-218.


\end{thebibliography}
\end{document}